\documentclass{amsart}
\usepackage{amssymb, latexsym}
\pdfoutput=1
\usepackage{graphicx}
\usepackage{listings}
\begin{document}
\title{Beta Distribution of  Long Memory Sequences}
\author{Robert Kimberk }
\curraddr{ Smithsonian Astrophysical Observatory\\  Cambridge Ma. \\ U.S.A.}
\email{rkimberk@cfa.harvard.edu}
\maketitle

\begin{abstract}
Three long memory models,  ARFIMA, Timmer and Konig 1995, and a circular convolution model based on Wold's representation theorem are examined. Each model is shown to produce sequences with nonstationary generalized beta  marginal distributions. It is demonstrated that  the variance  divided by the squared range of the sequence is stationary, and is a function of the shape parameter of the resulting symmetric beta distribution. Using the Wold model, a simple  matrix distribution transformation is given that maps the normal components of the long memory model onto the beta distribution.
\end{abstract}

\section{Introduction}

Two properties of long memory sequences, produced by the three models, obscure the probability  distribution: serial correlation, and nonstationarity of the variance and mean.  Serial correlation of the sequence impedes sampling of the full extent of the distribution. Generation of multiple subsequences derived from different random seeds  allow sampling of most of the possible values of the model. Min-max normalization of the same subsequences transforms the distribution to the unit interval $\lbrack 0,1 \rbrack$, mitigating nonstationarity, and preserving the distribution type. Concatenation of the processed subsequences facilitates the generation a complete empirical density. These operations allow a representation of the sequence distribution type as a symmetric beta distribution. It is inferred from the resulting distribution that the preprocessed distribution is a nonstationary generalized beta distribution.

The ratio of  variance and squared range of  long memory sequences  produced by the  three models is shown to be stationary. As well, the shape parameter of the symmetric  beta distribution is a function of the ratio. 

The three models tested are Granger and Joyeux's ARFIMA model, Timmer and Konig's 1995 model, and a model based on Herman Wold's moving average model. Kernal density plots  are derived from the ARFIMA model for truncated normal, Wigner semi-circle, uniform, and arcsine densities. These densities are  identifiable members of the symmetric beta distribution family.
The four beta  densities produced by the Timmer and Konig model and the Wold model are statistically identical and are not given in a figure.

A model generating long memory sequences based on Herman Wold's convolution theorem is given. The model consists of the matrix product of a circulant matrix derived from the spectral trend of the desired long memory sequence, and a vector of zero mean, unit variance, normal random samples. The simplicity of the linear algebra  model facilitates understanding of the distribution transformation of the normal distribution onto the beta distribution using the inverse circulant matrix.

Code written in the R language is used as the exploratory component of this paper, and examples of the code are provided. 
\section{The Shape Parameter of the Symmetric Generalized Beta Distribution}

The variance of a four parameter generalized beta distributed sequence is defined as
\begin {equation}
\sigma^2 = (b-a)^2 \frac {\alpha \beta}{ (\alpha + \beta)^2(\alpha +\beta +1)}
\end{equation} 
where $(b-a)^2 $ is the  squared range of the sequence, and $(\alpha,\beta)$ are the two beta distribution shape parameters. As all the beta distributions described in this paper are symmetric, which occurs when $\alpha$ and $\beta$ are equal, this paper combines the two shape parameters into the $\alpha$ symbol such that $\alpha + \beta = 2\alpha$ and $\alpha \beta = \alpha^2$ .  Merging the shape parameters in equation (1), dividing both sides by the squared range, and simplifying the resulting equation gives the following
\begin{equation}
\frac{\sigma^2 }{(b-a)^2} = \frac{1}{8\alpha +4}
\end{equation}
The shape parameter, $\alpha$, of the four beta distributions are as follows,  for arcsine distribution $\alpha= \frac{1}{2}$, 
uniform distribution $\alpha = 1$, Wigner semi circle distribution $\alpha = \frac{3}{2}$ and the truncated normal distribution $\alpha > 5$.

A similar  equation  is Tiberiu Popoviciu's variance inequality on bounded probability distributions \cite{tP35}
\begin{equation}
\frac{\sigma^2}{(b-a)^2} \leqq \frac{1}{4}
\end{equation}

All sequences in this paper are of bounded support and  elements of the sequences are finite. Popoviciu's theorem constrains $\alpha$ in equation (2) such that $0 \leqq \alpha < \infty$.

 A property of 
\begin{equation}
 \frac{\sigma^2}{(b-a)^2}
\end{equation}
is scale invariance. A homogeneous function $f$ of degree $ p$  is a function of a sequence $ x $ and scalar $t$ such that \cite{tA65}

\begin {equation}
f(t(x)) = t^p f(x)
\end{equation}
 It is easy to see that  variance, and squared range are both homogeneous functions of degree two. The ratio of variance to squared range remains constant when sequence $x$ is scaled by scalar $t$.

\section{Min-Max Normalization}

Min-max normalization is a scale and  location transformation that preserves the  type of distribution family \cite{aK37, bG54} , and maps the distribution to the unit interval $\lbrack 0,1 \rbrack$. Subtract the minimum value of sequence, $ min(x)$, from each sequence element, and divide each element of the result by the range, $range(x)$, to produce the min max normalized sequence.

\begin{equation}
\frac{x- min(x)}{range(x)}
\end{equation}

In terms of the sequence $x$, with  range  $\lbrack a,b\rbrack$, and with $ a = min(x)$, $b = max(x)$, min max normalization  is described as

\begin{equation}
\frac{\lbrack a,b\rbrack-a}{b-a} = \frac{\lbrack 0,b-a \rbrack}{b-a} = \lbrack 0,1 \rbrack
\end{equation}

For symmetric beta distributions with support on the unit interval $\lbrack 0,1 \rbrack$ the mean is one half, and the ratio of variance to squared range is equal to  the variance. The variance of the normalized four test beta distributions are: arcsine distribution $var = \frac{1}{8}$, uniform distribution $var = \frac{1}{12}$, Wigner semicircle distribution $var = \frac{1}{16}$, and truncated normal distribution var = small.

 \section{Power Spectral Density Trend of a Long Memory Sequence}

The power spectral density of a long memory sequence is characterized by random values distributed about a deterministic trend that is proportional to the absolute value of the Fourier frequency variable $f$ raised to a power $-\beta$, and scaled by a constant $k$

\begin{equation}
k|f|^{-\beta}  \   \  \ where \  \ 0<k,\  \ and \ \ 0 \leqq \beta 
\end{equation}
 
The power spectral density trend is a power law function producing a linear graph with a slope of $-\beta$, when plotted with logarithmic coordinate axes.

While $\beta$ is a convenient descriptor for a long memory sequence it is not sufficient to describe the probability  distribution of the long memory sequence. In general range is not defined by  $\beta$, with the exception of a random walk, where the range is defined by  Khinchin's  Law of the Iterated Logarithm. However there is a loose relation between the power spectral density $\beta$ and the beta distribution shape parameter $\alpha$. As $\beta$ increases $\alpha$ decreases.

\begin{figure}[b]
\centering\includegraphics[scale=.62]{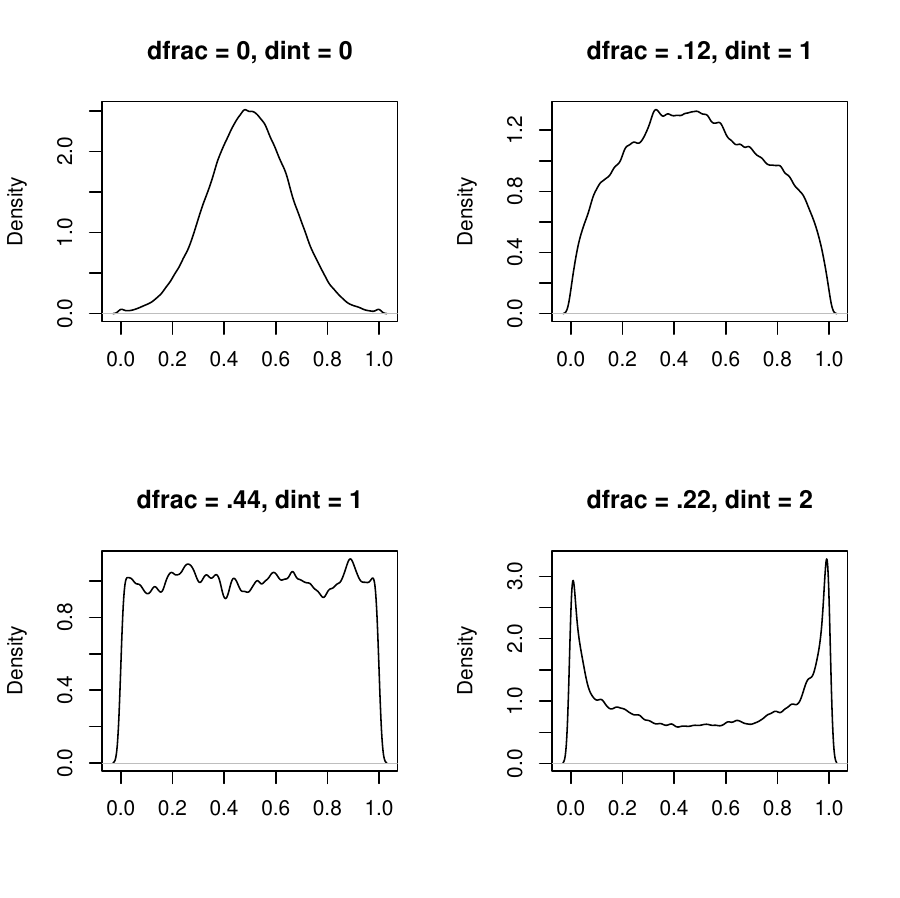}
\caption{Beta Densities from ARFIMA.SIM}
\label{Fi:001}
\end{figure}

\section{ Three Long Memory Models}

The four beta densities presented in figure 1 are representative of the densities generated by each of the three models. In each model a sequence of normally distributed values is transformed  into a  long memory sequence with a beta type distribution.The R code in the following four programs may be copied  and pasted into a R editor, and after the indicated parameters are changed in the code, the program may be  run. 

\subsection{ARFIMA}
The densities in figure 1 are derived from a fractionally differenced ARFIMA(0,d,0) model, where the d = dint + dfrac  parameters are given above each graph of density. Description of the model may be found in the following references \cite{cG80,jH81}.

The R program to generate the densities is as follows. The arfima.sim function parameters dint and dfrac may be changed to the values shown in figure 1 to derive the four beta type distributions. 
\vspace{0.60in}

\begin{verbatim}

#arfima.sim produces beta distribution
#load arfima package
n = 1e3
t = 0
i = 1
for ( i in 1:200 ) {
set.seed( i* 7)
q = arfima.sim( n, model = list(dfrac = .22, dint = 2, 
sigma2 = 1 ),  rand.gen = rnorm )
z = ( q - min( q ) ) / (max(q) -  min( q ) )
t = c( t, z ) }
var(t)
d = density( t, bw = .01 )
plot( d, lwd = 3 )

\end{verbatim}

Note that the distribution of the innovations is set to normal by rand.gen = rnorm, and the densities ploted are kernal densities.

\subsection{TK95}
A power law noise model created by Timmer and Konig in 1995. The model, based on function (8) of this paper,  is  fully explained in the paper by Timmer and Konig \cite{TK95}. 

The R program to generate the four densities as shown in figure 1 is as follows. The alpha argument of the TK95 function may be changed to 0 for a truncated normal density, 2.2 for a Wigner semi circle density, 2.9 for an uniform density, and 8 for an arcsine density.
\begin{verbatim}
#TK95 produces beta distribution
# load package RobPer
n = 1e3
i = 1
y = 0
for (i  in  1:200) {
set.seed(i*7)
x = TK95(n,alpha = 8)
x = (x - min(x)) / (max(x) - min(x))
y = c(y,x) }
var(y)
d = density(y, bw = .01)
plot(d, lwd  = 3)
\end{verbatim}
\subsection{Wold model}
A modification of Herman Wold's convolution model for a nonstationary stochastic process \cite{hW54}. It is a matrix product of a circulant matrix \cite{bB22}, derived from function (8) of this paper, and a random normal vector to a long memory  sequence. An advantage of this model is this, the inverse circulant matrix can be used to demonstrate the transformation of the normal distribution into a beta distribution.  A more  complete description of this model may be found in the paper by Kimberk, Carter, and Hunter \cite{RCH23}. The beta variable may may be changed to 0 for a truncated normal density, 2.17 for a Wigner semi circle  density, 2.87 for an uniform density, and 8.6 for an arcsine density.
\begin{verbatim}
#Wold's model produces beta distribution
#load package magic 
beta = 2.87
n = 1e3
b  =  0
i  =  1
y =  0
for (i in 1:100) {
set.seed(200 + i*7)
epsilon = rnorm(n+1, mean = 0, sd = 1)
freq = c(seq(-1/2, -1/n, by = 1/n), 1e-5, 
seq(1/n, 1/2, by = 1/n))
density = (abs(freq))^(-beta/2)
b = Mod(((n^(-1/2))* fft(density, inverse = TRUE )))
C = circulant(b)
x = C %*% epsilon
x = (x - min(x)) / (max(x) - min(x))
y = c(y,x) }
var(y)
d = density(y, bw = .01)
plot(d, lwd = 3)

\end{verbatim}

\section{Distribution Transformation}
Let C be the square circulant matrix in the Wold model, and x be the column vector of a normally distributed instances of random variable X,  as  in Wold's model, section 5.3. Then y, the product of matrix C and vector x, will be a long memory sequence with a beta type distribution.
\begin{equation}
y = C  x
\end{equation}

If  $F_ B$  is the beta distribution function and $F_N$ is the normal distribution function then the following equation is the transformation of distribution functions \cite{kS23}.
\begin{equation}
F_B[y] = P(Y \leqq y)= P(C \epsilon \leqq y) = P(X \leqq C^{-1} y) = F_N[C^{-1}(y)]
\end{equation}
The four beta distribution functions of figure 2, corresponding to the four density functions of figure 1, are derived from the following R program implementing  equation (10). Note that the range of the truncated normal distribution (beta = 0) is shortened by changing  d  to d[450:550] in the last line of the following code.

\begin{figure}[h]
\centering\includegraphics[scale = .65]{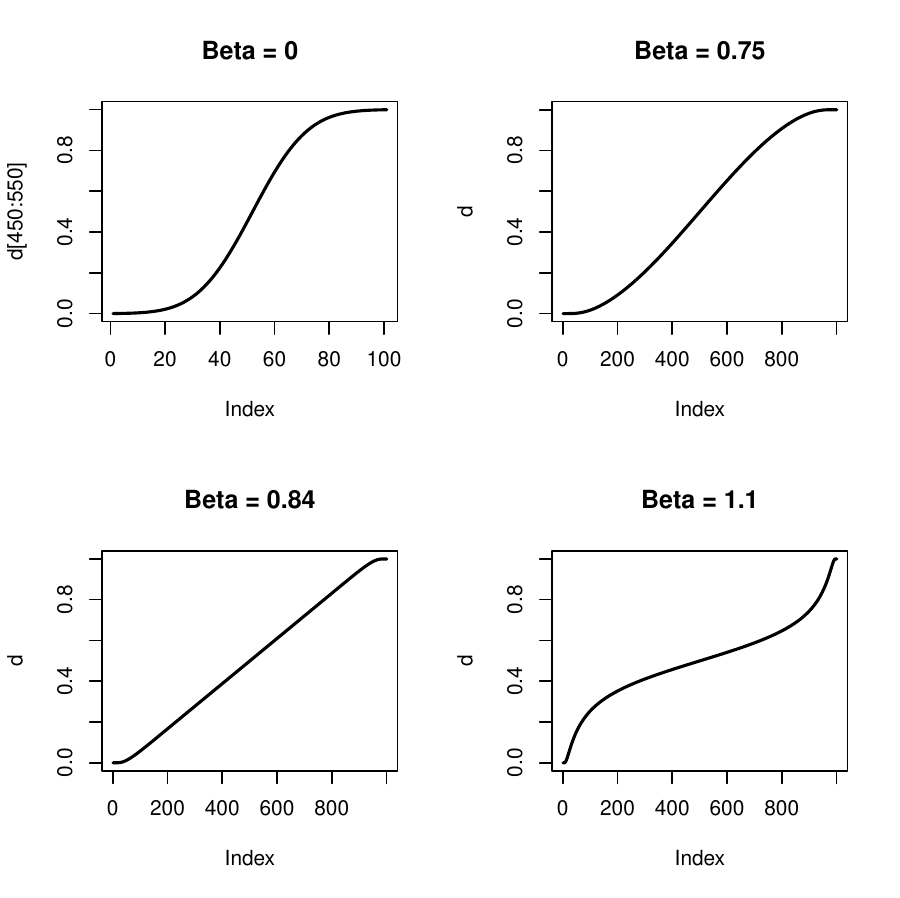}
\caption{Four Beta Distributions Functions of Transformed \ \  Normal Sequence}
\label{Fi:002}
\end{figure}
 \vspace{.4in}

\begin{verbatim}

# Distribution transformation
# load package magic
# beta: 0 = normal, 0.75 = Wigner semi circle, 
#0.84 = uniform, 1.1 = arcsine
beta = .75
n = 1000
freq = c(seq(-1/2, -1/n, by = 1/n), 1e-5, 
seq(1/n, 1/2, by = 1/n))
density = (abs(freq)) ^(-beta / 2)
b = Mod((n^(-1/2))*fft(density, inverse = TRUE))
C = circulant(b, doseq = TRUE)
y = seq(-1000, 1000, by = 2)
inverse_C = solve(C)
x = inverse_C %*%  y
d = pnorm(x, mean = 0, sd = 1)	
plot(d, type = "l", lwd = 2)

\end{verbatim}

\pdfoutput = 1

\vspace{.25in}

\end{document}